\documentclass[11pt,leqno,english]{article}
 \usepackage{latexsym}
 \usepackage{amssymb}
 \usepackage{amsmath}
 \usepackage{amsthm}
\usepackage{natbib}

 \usepackage[dvips]{graphics}
\usepackage{graphicx}

 \newcommand{ \un }{\mathbb{I} }
 \newcommand{ \p }{\mathbb{P} }
 \newcommand{ \pa }{\mathbb{P}^{\alpha} }
  \newcommand{ \pao }{\mathbb{P}^{\alpha}_{0} }
\newcommand{ \pam }{\mathbb{P}^{\alpha}_{m_n}}
 \newcommand{ \E }{\mathbb{E}}
 \newcommand{ \Ea }{\mathbb{E}^{\alpha}}
 \newcommand{ \Eam }{\mathbb{E}^{\alpha}_{m_n}}

\newcommand{ \F }{ \mathbb{F} }
 
 \newcommand{ \Z }{ \mathbb{Z} }
 \newcommand{\N}{ \mathbb{N} }
 
 \newcommand{ \V }{\textrm{Var} }

 \newcommand{ \Ct }{ \mathcal{C} }
 
 \newcommand{ \f }{ \mathcal{F} }

 \newcommand{ \tm }{ m }
 \newcommand{ \tM }{ M }
\newcommand{ \tmo }{ m_n }

\newcommand{ \tc }{ \tilde{c} }
\newcommand{ \Ie }{ I_{\eta_0} }
\newcommand{ \eto }{ \eta_0 }
\newcommand{ \tcu }{\tc_{3}(\beta) }
\newcommand{ \no }{\textrm{o} }
\newcommand{ \teb }{\theta_{\beta}(x) }

\newcommand{ \A}{ \mathcal{A} }
 \newcommand{ \B}{ \mathcal{B} }

 \newcommand{ \lo }{ \mathcal{L} }

\newtheorem{The}{{\bf Theorem}}[section]
\theoremstyle{definition}
\newtheorem{Def}[The]{{\bf Definition}}
\theoremstyle{plain}
 \newtheorem{Lem}[The]{Lemma}
 
 \newtheorem{Pro}[The]{\bf Proposition}
 \theoremstyle{definition}
\newtheorem{Rem}[The]{{\bf Remark}}

 \newenvironment{Pre}{\noindent \textbf{Proof.} \\ }{$\
 \blacksquare$}

\makeatletter

\@addtoreset{equation}{section} \makeatother

\title{Almost sure estimates for the concentration neighborhood of Sinai's walk
\\ \vspace{1cm}
 \large{Pierre Andreoletti} $^{\dag}$,\footnote{ Laboratoire Analyse-Topologie-Probabilit\'es - C.N.R.S. UMR 6632  
Universit\'e Aix-Marseille I, 
(Marseille France). \newline \vspace{0.1cm}  $\quad$  MSC 2000 60G50; 60J55. \newline \vspace{0.5cm} \textit{Key words and phrases :  Random environment, random walk,
Sinai's regime, local time, concentration.} }   }

\begin{document}

\bibliographystyle{unsrtnat}
\maketitle

\noindent  $^{\dag}$ Universit\'e Aix-Marseille I, Centre de math\'ematiques et d'informatique, 39 rue F. Joliot-Curie, 
13453 Marseille cedex 13, 
France. e-mail : \texttt{andreole@cmi.univ-mrs.fr}

\noindent \\ \textbf{Abstract:} We consider Sinai's random walk in
random environment. We prove that infinitely often (i.o.) the size of the concentration neighborhood of this random walk is almost surely bounded. As an application we get that i.o. the maximal distance between two favorite sites is almost surely bounded.

\section{Introduction and results}

In this paper we are interested in Sinai's walk i.e a one
dimensional random walk in random environment with three
conditions on the random environment: two necessaries hypothesis
to get a recurrent process (see \cite{Solomon}) which is not a
simple random walk and an hypothesis of regularity which allows us
to have a good control on the fluctuations of the random
environment. 
The asymptotic behavior of such walk was discovered by
\cite{Sinai} : this walk is sub-diffusive and at an instant $n$ it is localized in the neighborhood of a well
defined point of the lattice. The correct almost sure behavior of this walk, originally studied by \cite{Deh&Revesz}, have been checked by the remarkable precise results of \cite{HuShi2}. 
We denote Sinai's walk $(X_{n}, n \in \N)$, let us define the local time $\lo$, at $k$ $(k\in \Z)$ within the
interval of time $[1,T]$ ($T \in \N^*$) of $(X_n,n\in \N)$
\begin{eqnarray}
\lo\left(k,T\right) \equiv \sum_{i=1}^T \un_{\{X_i=k\}} .
\end{eqnarray}
$\un$ is the indicator function ($k$ and $T$ can be deterministic or random variables). Let $V\subset \Z$, we denote
\begin{eqnarray}
\lo\left(V,T\right) \equiv \sum_{j \in V} \lo\left(j,T\right)
=\sum_{i=1}^T\sum_{j \in V} \un_{\{X_i=j\}} .
\end{eqnarray}
\noindent Now, let us introduce the following random variables
\begin{eqnarray}
& & \lo^*(n)=\max_{k\in \Z}\left(\lo(k,n)\right) ,\
\F_{n}=\left\{k \in \Z,\ \lo(k,n)=\lo^*(n) \right\} , \\
&& Y_n=\inf_{x \in \Z}\min\left\{k>0\ :\ \lo([x-k,x+k],n) \geq n/2\right\}.
\end{eqnarray}
 $\lo^*(n)$ is the maximum of the local times (for a given instant
$n$), $\F_{n}$ is the set of all the favourite sites and $Y_{n}$ is the size of the interval where the walk spends more than a half of its time.
The first almost sure results on the local time are given by \cite{Revesz}, he notices and shows in a special case that $\lo^*$ can be very big (see also \cite{Revesz1}), then \cite{Shi} proves the result in the general case (we recall this result here : Theorem \ref{th0}). 
About $\F_{n}$, in \cite{HuShi0} it is proven, that the maximal favorite site is almost surely transient and that it has the same almost sure behavior as the walk itself (see also \cite{Shi1}).
Until now, the random variable $Y_{n}$ has not been studied a lot for Sinai's walk.  
In \cite{Pierre2} it is proven, that in probability, this random variable is very small comparing to the typical fluctuations of Sinai's walk. Here we are interested in the almost sure behavior of $Y_{n}$. We prove that the ''liminf'' of this random variable is almost surely bounded. We will see that the result we give for $Y_{n}$ implies the result of R\'ev\'esz about $\lo^*$ and have interesting consequence on the favorite sites.

\subsection{Definition of Sinai's walk}

Let $\alpha =(\alpha_i,i\in \Z)$ be a sequence of i.i.d. random
variables taking values in $(0,1)$ defined on the probability
space $(\Omega_1,\f_1,Q)$, this sequence will be called random
environment. A random walk in random environment (denoted
R.W.R.E.) $(X_n,n
\in\N)$ is a sequence of random variable taking value in $\Z$, defined on $( \Omega,\f,\p)$ such that \\
$ \bullet $ for every fixed environment $\alpha$, $(X_n,n\in \N)$
 is a Markov chain with the following transition probabilities, for
 all $n\geq 1$ and $i\in \Z$
 \begin{eqnarray}
& & \p^{\alpha}\left[X_n=i+1|X_{n-1}=i\right]=\alpha_i, \label{mt} \\
& & \p^{\alpha}\left[X_n=i-1|X_{n-1}=i\right]=1-\alpha_i \equiv
\beta_i. \nonumber
\end{eqnarray}
We denote $(\Omega_2,\f_2,\pa)$ the probability space
associated to this Markov chain. \\
 $\bullet$ $\Omega = \Omega_1 \times \Omega_2$, $\forall A_1 \in \f_1$ and $\forall A_2 \in \f_2$,
$\p\left[A_1\times
A_2\right]=\int_{A_1}Q(dw_1)\int_{A_2}\p^{\alpha(w_1)}(dw_2)$.

\noindent \\ The probability measure $\pa\left[\left.
.\right|X_0=a \right]$  will be  denoted $\pa_a\left[.\right]$,
 the expectation associated to $\pa_a$: $\Ea_a$, and the expectation associated to $Q$:
 $\E_Q$.

\noindent \\ Now we introduce the hypothesis we will use in all
this work. The two following hypothesis are the necessaries
hypothesis
\begin{eqnarray}
 \E_Q\left[ \log
\frac{1-\alpha_0}{\alpha_0}\right]=0 , \label{hyp1bb} \label{hyp1}
\end{eqnarray}
\begin{eqnarray}
\V_Q\left[ \log \frac{1-\alpha_0}{\alpha_0}\right]\equiv \sigma^2
>0 . \label{hyp0}
\end{eqnarray}
 \cite{Solomon} shows that under \ref{hyp1} the
process $(X_n,n\in \N)$ is $\p$ almost surely recurrent and
\ref{hyp0} implies that the model is not reduced to the simple
random walk. In addition to \ref{hyp1} and \ref{hyp0} we will
consider the following hypothesis of regularity, there exists $0<
\eta_0 < 1/2$ such that
\begin{eqnarray}
& & \sup \left\{x,\ Q\left[\alpha_0 \geq x \right]=1\right\}= \sup
\left\{x,\ Q\left[\alpha_0 \leq 1-x \right]=1\right\} \geq \eta_0.
\label{hyp4}
\end{eqnarray}

\noindent We call \textit{Sinai's random walk} the random walk in
random environment previously defined with the three hypothesis
\ref{hyp1}, \ref{hyp0} and \ref{hyp4}.

 \subsection{Main results}

\begin{The} \label{th1}
Assume \ref{hyp1bb}, \ref{hyp0} and \ref{hyp4} hold,
  there exists $c_{1}\equiv c_{1}(Q)>0$ such that
\begin{eqnarray}
& & \p\left[ \liminf_{n} Y_{n} \leq c_{1}
 \right]  =1. 
\end{eqnarray}
\end{The}
\noindent This first result prove that, almost surely, one can find a subsequence such that the size of the neighborhood where the walk spend more than a half of its time is bounded from above by a constant depending only on the distribution of the random environment. 

\noindent As a corollary  we get the following result originally due to \cite{Revesz} but which proof have been performed in the general case by \cite{Shi} :
\begin{The}   \label{th0}
 Assume \ref{hyp1bb}, \ref{hyp0} and \ref{hyp4} hold,
  there exists $c_{2}\equiv c_{2}(Q)>0$ such that
\begin{eqnarray}
& & \p\left[ \limsup_{n} \frac{\lo^*(n)}{n} \geq c_{2}
 \right]  =1.
\end{eqnarray}
\end{The}

\noindent In \cite{Pierre2} we were also interested in the size of the interval, centered on the point of localisation defined by \cite{Sinai}, where the walk spends  an arbitrary proportion of time $n$ (see, for example, Theorem 3.1 in \cite{Pierre2}). It is proven that the size of this intervall is once again negligible comparing to the typical fluctuation of the walk. Here we are interested  in the following random variable, let $0 \leq \beta<1$
\begin{eqnarray}
Y_{n,\beta}=\inf_{x \in \Z}\min\left\{k>0\ :\ \lo([x-k,x+k],n) \geq \beta n \right\},
\end{eqnarray}
notice that $Y_{n}\equiv Y_{n,1/2}$, we get the following result

\begin{The}  \label{th2}
Assume \ref{hyp1bb}, \ref{hyp0} and \ref{hyp4} hold,
  there exists $c_{3}\equiv c_{3}(Q)>0$ such that for all $0\leq \beta<1$
\begin{eqnarray}
& & \p\left[ \liminf_{n} Y_{n,\beta} \leq c_{3}(1-\beta)^{-2}
 \right]  =1. 
\end{eqnarray}
\end{The}
\noindent
We notice that, when $\beta$ get close to one, meaning that we look for the size of an interval where the local time is close to $n$, the size of this interval grows like $1/(1-\beta)^{2}$. Of course this result implies Theorem \ref{th1}. We will explain in detail this $1/(1-\beta)^{2}$ dependence.\\
\noindent As an application we get the following result about the maximal distance between two favorite sites, 
\begin{The}  \label{th3}
Assume \ref{hyp1bb}, \ref{hyp0} and \ref{hyp4} hold,
  there exists $c_{4}\equiv c_{4}(Q)>0$ such that 
\begin{eqnarray}
& & \p\left[ \liminf_{n} \max_{(x,y)\in \F_{n}^2} |x-y| \leq c_{4}
 \right]  =1. 
\end{eqnarray}
\end{The}
\noindent

\noindent We get that infinitely often the maximal distance between two favorite sites is almost surely bounded, notice that this implies also that, almost surely, there is only a finite number of favorite sites at step $n$ infinitely often.

\subsection{About the proof of the results}

We have used a similar method of \cite{Pierre2}, and also an extension for Sinai's walk of Propositon 3.1 of \cite{GanShi}. We will give the details of proof in such a way  the reader understand the $(1-\beta)^{-2}$ dependance occurring in Theorem \ref{th2}. However some details of proof, already present in \cite{Pierre2}, have not been repeated here.

\noindent \\ This paper is organized as follows. In section 2 we
give the proof of Theorems \ref{th1} to \ref{th2}, in section 3 we prove Theorem \ref{th3}, finally in section 4 we point out remarks and open questions.
In the appendix we give the needed estimate for the environment and detail the proof of some of it.

\section{Proof of Theorems \ref{th1}-\ref{th2}}

We have point out that Theorem \ref{th2} implies the two other (\ref{th1} and \ref{th0}), so the main part of this section is to prove this Theorem. Notice, that Theorem \ref{th1} is a special case of Theorem \ref{th2} taking $\beta=1/2$, at the end of the section we will  explain why we also get Theorem \ref{th0}.

\noindent \\ To prove Theorem \ref{th2} we begin with  the following elementary remark  :
By definition we have 
\begin{eqnarray}
\liminf_{n} Y_{n,\beta} \leq c_{3}(1-\beta)^{-2} \iff \bigcap_{N} \bigcup _{n \geq N}\left\{ Y_{n,\beta} \leq c_{3}(1-\beta)^{-2}\right\},
\end{eqnarray}
denote $\tcu \equiv  c_{3}(1-\beta)^{-2} $ and $\teb=[x-\tcu,x+\tcu]$, we have the inclusion 
\begin{eqnarray}
\left \{ \max_{x}\lo\left(\teb,n\right) \geq \beta n \right\}  \subseteq \left\{ Y_{n,\beta} \leq \tcu \right \},
\end{eqnarray}
so we get that 
\begin{eqnarray}
\p\left[\liminf_{n} Y_{n,\beta} \leq \tcu \right] &\geq&  \p\left[ \bigcap_{N} \bigcup _{n \geq N} \left \{  \max_{x}\lo\left(\teb,n\right) \geq \beta n \right\} \right] \nonumber \\
&\equiv &  \p\left[ \limsup_{n} \frac{ \max_{x}\lo\left(\teb,n\right)}{n} \geq \beta  \right].
\end{eqnarray}
\noindent 
To get the result it is enough to prove the two following Propositions :
\begin{Pro} \label{Prop1} Let $(\phi(n),n)$ be a strictly positive sequence such that $\lim_{n \rightarrow \infty} \phi (n)=+\infty$, for all $0 \leq \beta<1$ we have
\begin{eqnarray}
\p\left[ \limsup_{n} \frac{ \max_{x}\lo\left(\teb,n\right)}{\phi(n)}= \textrm{const} \in [0,\infty]  \right]=1.
\end{eqnarray}
\end{Pro}
\noindent and 
\begin{Pro} \label{Prop21}   For all $0 \leq  \beta<1$ we have
\begin{eqnarray}
\p\left[  \frac{ \max_{x}\lo\left(\teb,n\right)}{n} \geq \beta   \right] >0 .
\end{eqnarray}
\end{Pro}
\noindent
Notice that Proposion \ref{Prop1} is a simple extension  for Sinai's walk of Proposition 3.1 of \cite{GanShi}, as one can find the details of the proof in the referenced paper, we just explain why it works in our case :
\\
\subsection{Proof of Proposition \ref{Prop21}}
Define $f(\alpha,(X_{m}))=\limsup_{n} \frac{ \max_{x}\lo\left(\teb,n\right)}{\phi(n)}$, following the method of \cite{GanShi} it is enough to prove the two following facts : \textit{Fact} 1 for $Q$-a.a. $\alpha$ $f(\alpha,(X_{m}))$ is constant for $\pa$-a.a. realizations of $(X_{n},n)$ and \textit{Fact} 2 $f(\alpha)\equiv f(\alpha,(X_{m}))$ is a constant for $Q$-a.a. $\alpha$.
The key point for the proof of this two facts is that for all $x\in\Z$ ($T_{x}<+\infty$  $\pa$-a.s for $Q$-a.a. $\alpha$)  because Sinai's walk is $\p$-a.s recurrent. So we can apply the three steps of the proof of \cite{GanShi} (pages 168-169) : the two first provide Fact 1, the third one Fact 2. 
Notice that here we need a result for $\max_{x}\lo\left(\teb,n\right)$, with $\teb$ a finite  interval, whereas in \cite{GanShi} $\max_{x}\lo\left(x,n\right)$ is studied, however this difference does not change the computations. 
$\blacksquare$

\subsection {Proof of Proposition \ref{Prop2}}
To prove this Proposition we use a quite similar method of \cite{Pierre2}, first let us recall  the following decomposition of the measure $\p$, let $\Ct_{n} \in \sigma \left(X_{i}, i \leq n \right) $ and $G_{n} \subset \Omega_{1}$, we have :
\begin{eqnarray}
\p\left[\Ct_{n}\right] &\equiv & \int_{\Omega_{1}}Q(d\omega)\int_{\Ct_{n}} d\p^{\alpha(\omega)} \\
 & \geq &  \int_{G_{n}}Q(d\omega)\int_{\Ct_{n}} d\p^{\alpha(\omega)}.
\end{eqnarray}
So assume that for all $\omega \in G_{n}$ and $n$, $\int_{\Ct_{n}} d\p^{\alpha(\omega)} \equiv d_{1}(\omega,n)>0$ and assume that $Q[G_{n}] \equiv d_{2}(n)>0$ we get that for all $n$
\begin{eqnarray}
\p\left[\Ct_{n}\right] \geq d_{2}(n) \times \min_{w\in G_{n}}(d_{1}(w,n))>0.
\end{eqnarray}
So choosing $\Ct_{n}= \left\{\max_{x}\lo\left(\teb,n\right) \geq \beta n \right\}$,  we have to extract from $\Omega_{1}$ a subset $G_{n}$ sufficiently small to get that $\min_{w\in G_{n}}(d_{1}(w,n))>0$ (Proposition \ref{Prop2}) but sufficiently large to have $d_{2}(n)>0$ (Proposition \ref{profondab}) .
The largest part of the proof is to construct such a $G_{n}$ (Section \ref{sec11} and Appendix B).

\subsubsection{ Construction of $G_{n}$ (arguments for the random environment) \label{sec11}}

For completeness we begin with some basic notions originally introduced by \cite{Sinai}.
\noindent \\ \textbf{The random potential and the valleys}
\noindent \\ Let
\begin{eqnarray}
\epsilon_i \equiv \log \frac{1-\alpha_i}{\alpha_i},\ i\in \Z,
\end{eqnarray}
define :
\begin{Def} \label{defpot2} The random potential $(S_m,\  m \in
\Z)$ associated to the random environment $\alpha$ is defined in the following way: for all $k$ and $j$, if $k>j$
\begin{eqnarray} 
 && S_{k}-S_{j}=\left\{ \begin{array}{ll} \sum_{j+1\leq i \leq k} \epsilon_i, &  k\neq0, \\
  - \sum_{j \leq i \leq -1} \epsilon_i , &  k=0  , \end{array} \right. \nonumber \\
&& S_{0}=0, \nonumber
 \end{eqnarray}
and symmetrically if $k<j$.
\end{Def}

\begin{figure}[h]
\begin{center}
\input{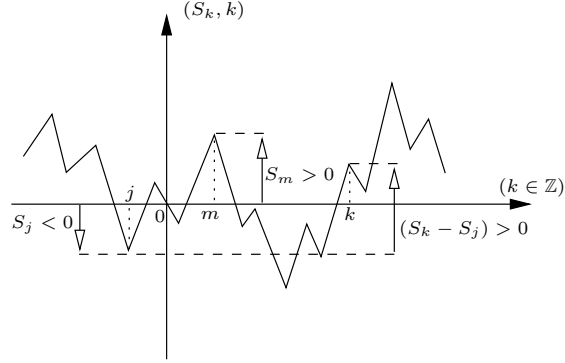} \caption{Trajectory of the random potential} \label{fig5}
\end{center}
\end{figure}

\noindent
\begin{Rem}  
using Definition \ref{defpot2} we have :
\begin{eqnarray} 
 S_k=\left\{ \begin{array}{ll} \sum_{1\leq i \leq k} \epsilon_i, &  k=1,2,\cdots , \\
   \sum_{k \leq i \leq -1} \epsilon_i , &  k=-1,-2,\cdots  , \end{array}  \right. \label{defmerde}
 \end{eqnarray}
however, if we use \ref{defmerde} for the definition of $(S_{k},k)$, $\epsilon_{0}$ does not appear in this definition and moreover it is not clear, when $j<0<k$, what the difference $S_k-S_{j}$ means (see figure \ref{fig5}).   
\end{Rem}

\begin{Def} \label{c2s2d1}
 We will say that the triplet $\{M',m,M''\}$ is a \textit{valley} if
 \begin{eqnarray}
& &  S_{M'}=\max_{M' \leq t \leq m} S_t ,  \\
& &  S_{M''}=\max_{m \leq t \leq
 \tilde{M''}}S_t ,\\
& & S_{m}=\min_{M' \leq t \leq M''}S_t \ \label{2eq58}.
 \end{eqnarray}
If $m$ is not unique we choose the one with the smallest absolute
value.
 \end{Def}

\begin{Def} \label{deprofvalb}
 We will call \textit{depth of the valley} $\{\tM',\tm,\tM''\}$ and we
 will denote it $d([M',M''])$ the quantity
\begin{eqnarray}
 \min(S_{M'}-S_{m},S_{M''}-S_{m})
 .
 \end{eqnarray}
 \end{Def}

\noindent  Now we define the operation of \textit{refinement}
 \begin{Def}
Let  $\{\tM',\tm,\tM''\}$ be a valley and let
  $\tM_1$ and $\tm_1$ be such that $\tm \leq \tM_1< \tm_1 \leq \tM''$
  and
 \begin{eqnarray}
 S_{\tM_1}-S_{\tm_1}=\max_{\tm \leq t' \leq t'' \leq
 \tM''}(S_{t'}-S_{t''}) .
 \end{eqnarray}
 We say that the couple $(\tm_1,\tM_1)$ is obtained by a \textit{right refinement} of $\{\tM',\tm,\tM''\}$. If the couple $(\tm_1,\tM_1)$ is not
 unique, we will take  the one such that $\tm_1$ and $\tM_1$ have the smallest  absolute value. In a similar way we
  define the \textit{left refinement} operation.
 \end{Def}

\begin{figure}[h]
\begin{center}
\input{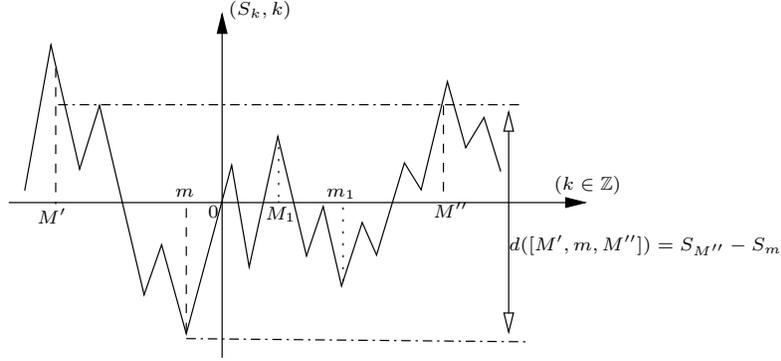} \caption{Depth of a valley and refinement operation} \label{fig4}
\end{center}
\end{figure}

\noindent \\  We denote $\log_2=\log \log $, in all this section we will suppose that $n$ is
large enough such that $\log_2 n$ is positive.

\begin{Def} \label{thdefval1b} Let $n>3$ and $\Gamma_n \equiv \log n+ 12 \log_2 n $, we say that a valley $\{\tM',\tm,\tM''\}$ contains $0$ and
is of depth larger than $\Gamma_n$ if and only if
\begin{enumerate}
\item  $ 0 \in [\tM',\tM'']$, \item $d\left([\tM',\tM'']\right)
\geq \Gamma_n $ , \item if $\tm<0,\ S_{\tM''}-\max_{ \tm \leq t
\leq 0
}\left(S_t\right) \geq 12 \log_2 n $ , \\
 if $\tm>0,\
S_{\tM'}-\max_{0 \leq t \leq
 \tm}\left(S_t\right) \geq 12 \log_2 n $ .
\end{enumerate}
\end{Def}

\noindent \textbf{The basic valley $\{{M_n}',\tmo,{M_n}\}$}

 We recall the notion of \textit{ basic valley } introduced
by Sinai and denoted here $\{{M_n}',\tmo,{M_n}\}$. The definition
we give is inspired by the work of \cite{Kesten2}. First let
$\{\tM',\tmo,\tM''\}$ be the smallest \textit{valley that contains
$0$ and of depth larger than} $\Gamma_n$. Here smallest means that
if we construct, with the operation of refinement, other valleys
in $\{\tM',\tmo,\tM''\}$ such valleys will not satisfy one of  the
properties of Definition \ref{thdefval1b}. ${M_n}'$ and ${M_n}$
are defined from $\tmo$ in the following way: if $\tmo>0$
\begin{eqnarray}
& & {M_n}'=\sup \left\{l\in \Z_-,\ l<\tmo,\ S_l-S_{\tmo}\geq
\Gamma_n,\
S_{l}-\max_{0 \leq k \leq \tmo}S_k \geq 12 \log_2 n \right\} ,\\
& & {M_n}=\inf \left\{l\in \Z_+,\ l>\tmo,\ S_l-S_{\tmo}\geq
\Gamma_n\right\} . \label{4.8}
\end{eqnarray}
if $\tmo<0$
\begin{eqnarray}
& & {M_n}'=\sup \left\{l\in \Z_-,\ l<\tmo,\ S_l-S_{\tmo}\geq
\Gamma_n\right\} , \\
& & {M_n}=\inf \left\{l\in \Z_+,\ l>\tmo,\ S_l-S_{\tmo}\geq
\Gamma_n,\ S_{l}-\max_{ \tmo \leq k \leq 0}S_k \geq 12 \log_2
n \right\} . \label{4.10}
\end{eqnarray}
if $\tmo=0$
\begin{eqnarray}
& & {M_n}'=\sup \left\{l\in \Z_-,\ l<0,\ S_l-S_{\tmo}\geq
\Gamma_n \right\} , \\
& &  {M_n}=\inf \left\{l\in \Z_+,\ l>0,\ S_l-S_{\tmo}\geq \Gamma_n
\right\} . \label{4.12}
\end{eqnarray}
\noindent  $\{{M_n}',\tmo,{M_n}\}$ exists
 with a $Q$ probability as close to one as we need. In fact it is not
 difficult to prove the following lemma

\begin{figure}[h]
\begin{center}
\input{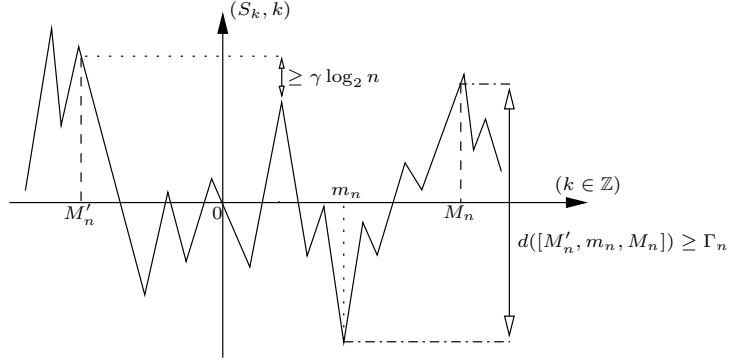} \caption{Basic valley, case $m_{n}>0$ } \label{thfig4}
\end{center}
\end{figure}

\begin{Lem} \label{moexiste} Assume  \ref{hyp1}, \ref{hyp0} and \ref{hyp4} hold, for all $n$ we have
\begin{eqnarray}
Q\left[\{{M_n}',\tmo,{M_n}\} \neq \varnothing \right] =
1-\no(1).
\end{eqnarray}
we denote $\no(1)$ a positive function of $n$ such that $\lim_{n \rightarrow  \infty } \no(1)=0$.
\end{Lem}

\begin{Pre}
One can find the proof of this Lemma in Section 5.2 of \cite{Pierre2}.
\end{Pre}

\noindent Let $x \in \Z$, define
\begin{eqnarray}
& &  T_{x}=\left\{\begin{array}{l} \inf\{k\in\N^*,\ X_k=x \}
\\ + \infty \textrm{, if such }  k \textrm{ does not exist.}
 \end{array} \right. \label{3.7sec}
\end{eqnarray}

\begin{Def} \label{superb}  Let $c_{0}>0$, $c_{3}>0$, $0 \leq \beta<1$ and $\omega \in \Omega_1$, we will say that $\alpha \equiv
\alpha(\omega)$ is a \textit{good environment} if there exists
$n_0$ such that for all $n\geq n_0$ the sequence $(\alpha_i,\ i
\in \Z)=(\alpha_i(\omega),\ i \in \Z)$ satisfies
 the properties \ref{3eq325} to \ref{8eq36}
\begin{eqnarray}
& \bullet &  \{{M_n}',\tmo, {M_n}\} \neq \varnothing,   \label{3eq325}\\
 & \bullet  & {M_n}'\geq ( \sigma^{-1} \log n)^2
 ,\ {M_n}\leq ( \sigma^{-1} \log n)^2  ,
\label{interminibb} \\
& \bullet & \Ea_{\tmo}\left[ \lo(\tilde{\Theta}(n,\beta),T_{\tmo})\right] \leq
 \frac{2 c_{0}}{\sqrt{\tcu}} . \label{8eq36}
\end{eqnarray}
where $\tilde{\Theta}(n,\beta)=[M_{n}',M_{n}'+1,\cdots,m_{n}-\tcu]\cup[m_{n}+\tcu,m_{n}+\tcu+1,\cdots,M_{n}]$, recall that $\tcu=c_{3}(1-\beta)^{-2}$.
\end{Def}


\noindent  Define the \textit{set of good environments}
\begin{eqnarray}
G_n \equiv G_{n}(c_{0},c_{3},\beta)=\left\{\omega \in \Omega_1,\ \alpha(\omega) \textrm{ is a }
\textit{ good environment} \right\}.
\end{eqnarray}
$G_n$ depends on $c_0$ and $n$, however  we do not
make explicit its $c_{0}$ and $\beta$ dependence.

\begin{Pro} \label{profondab} There exists a numerical constant $d_{2}>0$ such that if   \ref{hyp1bb}, \ref{hyp0} and \ref{hyp4} hold,  there exists $c_{0}\equiv c_{0}(Q)>0$ and  $n_0$ 
such that for all $c_{3}>0$, $0 \leq \beta<1$ and $n>n_0$
\begin{eqnarray}
Q\left[ G_n\right]  \geq 1/2 .
\end{eqnarray}
\end{Pro}

\begin{Pre}
We have to prove that the three properties \ref{3eq325}-\ref{8eq36} are true with a $Q$ probability strictly positive. In \cite{Pierre2} we prove that the two first properties are true with a probability close to one, so we only have to prove that the third one is true with a probability larger than $1/2$, this is done in the Appendix B.
\end{Pre}

\subsubsection{Argument for the walk (environment fixed $\alpha \in G_{n}$)}

In this section we assume that $n$ is sufficiently large such that Proposition \ref{profondab} is true and we assume also that the random environment is fixed and belongs to $G_{n}$ (denoted $\alpha \in G_{n}$). 

\begin{Pro} \label{Prop2}   For all $0 \leq \beta<1$, $n$ large enough and $\alpha \in G_{n}$ we have
\begin{eqnarray}
\pao\left[   \max_{x}\lo\left(\teb,n\right) \geq \beta n  \right] >1/2,
\end{eqnarray}
recall that $\teb=[x-\tcu,x+\tcu]$, $c_{3}\equiv c_{3}(Q)>0$.
\end{Pro}

\begin{Pre}
\noindent To get this result, it is enough to prove that, 
\begin{eqnarray}
\pao\left[  \lo\left(\theta_{\beta}(\tmo),n\right) \geq \beta n   \right] >1/2,
\end{eqnarray}
we will prove the following equivalent fact 
\begin{eqnarray}
\pao\left[  \lo\left(\Theta(n,\beta),n\right)  \geq (1-\beta)n   \right] <1/2, \label{eqes1}
\end{eqnarray}
 where $\Theta(n,\beta)$ is the complementary of $\theta_{\beta}(\tmo)$ in $\Z$.

\noindent \\
First we recall the two following elementary results
\begin{Lem}\label{lem2b}  For all $n$ and $\alpha \in G_{n}$ we have
\begin{eqnarray}
&& \p^{\alpha}_0\left[\bigcup_{m=0}^{n}\left\{X_m\notin \left[M_n',M_{n}
\right]\right\}\right] =
\no(1).  \label{eded} \\
&& \p^{\alpha}_0\left[ T_{\tmo}> \frac{n}{(\log n)^4} \right]=
\no(1) \label{4.27}.
\end{eqnarray}
Recall that $\lim_{n \rightarrow  \infty } \no(1)=0$.
 \end{Lem}

\begin{Pre}
This is a basic result for Sinai's walk, it makes use Properties \ref{3eq325} and \ref{interminibb}. One can find the details of this proof in \cite{Pierre2} : Proposition 4.7 and Lemma 4.8. 
\end{Pre}

\noindent \\ First we use \ref{eded} to 
reduce the set $\Theta(n,\beta)$ to $\tilde{\Theta}(n,\beta)$ defined just after \ref{8eq36}, we get
\begin{eqnarray}
& & \pao\left[ \lo\left(\Theta(n,\beta),n\right) \geq (1-\beta) n
 \right]  \leq \pao\left[ \lo\left(\tilde{\Theta}(n,\beta),n\right) \geq
(1-\beta) n \right]+ \no(1)
\label{4eq8}.
\end{eqnarray}
Now using \ref{4.27} we get 
\begin{eqnarray}
\qquad \pao\left[\lo\left(\tilde{\Theta}(n,\beta),n\right) \geq (1-\beta) n
\right]  & \leq & \pao\left[\lo\left(\tilde{\Theta}(n,\beta),n\right) \geq
(1-\beta) n,\ T_{\tmo} \leq \frac{n}{(\log n)^4}  \right]+ \no(1).
\end{eqnarray}
Let us denote $N_0 \equiv \left[n (\log n)^{-4}\right] +1$ and
$1-\beta_n \equiv 1-\beta-N_0/n$. By the Markov property and the
homogeneity of the Markov chain we obtain 
\begin{eqnarray}
\qquad \pao\left[\lo\left(\tilde{\Theta}(n,\beta),n\right) \geq (1-\beta) n,\
T_{\tmo} \leq \frac{n}{(\log n)^4}  \right] & \leq &
 \pam\left[ \sum_{k=1}^{n}\un_{\left\{X_k \in
\tilde{\Theta}(n,\beta) \right\} } \geq (1-\beta_n) n
 \right] \label{4eq21} .
\end{eqnarray}
\noindent Let $j\geq 2 $, define the following return times
\begin{eqnarray*}
&& T_{\tmo,j } \equiv \left\{\begin{array}{l}  \inf\{k>T_{\tmo,j-1},\ X_k=\tmo \}, \\
 + \infty \textrm{, if such }  k \textrm{ does not exist.}
 \end{array} \right. \\
& & T_{\tmo,1}  \equiv T_{\tmo} \ (\textrm{see \ref{3.7sec}}).
\end{eqnarray*}
Since by definition $T_{ \tmo,n} > n $, 
$\left \{ \sum_{k=1}^{n}\un_{\left\{X_k \in
\tilde{\Theta}(n,\beta) \right\} } \geq   (1-\beta_n) n
 \right \}
 \subset  \left \{ \sum_{k=1}^{T_{
\tmo,n}}\un_{\left\{X_k \in \tilde{\Theta}(n,\beta)\right\} } \geq (1-\beta_n) n
\right \}$, then using the definition of the local time and the Markov inequality we get
\begin{eqnarray}
\pam\left[ \sum_{k=1}^{n}\un_{\left\{X_k \in
\tilde{\Theta}(n,\beta) \right\} } \geq
 (1-\beta_n) n \right] & \leq & 
 \pam\left[ \sum_{k=1}^{T_{\tmo},n}\un_{\left\{X_k \in
\tilde{\Theta}(n,\beta) \right\} } \geq
 (1-\beta_n) n \right]  \\
 & \leq &  \Eam\left[\lo\left(\tilde{\Theta}(n,\beta),T_{\tmo}\right)\right] (1-\beta_n)^{-1} ,
\end{eqnarray}
and we have used the fact that, by the strong Markov property, the
random variables $\lo\left(s,T_{ \tmo,i+1}-T_{\tmo,i}\right)$ $(0
\leq i \leq n-1)$ are $i.d.$.
Using  the property \ref{8eq36}, there exists $c_{0}$ such that
\begin{eqnarray}
\Eam\left[\lo\left(\tilde{\Theta}(n,\beta),T_{\tmo}\right)\right] \leq \frac{2 c_{0}(1-\beta)}{(c_{3})^{1/2}}.
 \label{4eq36}
\end{eqnarray}
 Collecting what we did
above, we finally get for $n$ sufficiently large 
\begin{eqnarray}
\pao\left[\lo\left(\Theta (n,\beta),n\right) \geq (1-\beta) n \right] \leq \frac{4 c_{0}}{(c_{3})^{1/2}}.
\end{eqnarray}
we get \ref{eqes1} choosing $c_{3}= 64( c_{0})^{2}$.
\end{Pre}
\noindent \\
This ends the proof of Theorem \ref{th2}, to get Theorem \ref{th0}, we remark that
\begin{eqnarray*}
\left\{Y_{n}\leq c_{1} \right\} & \equiv & \left \{ \inf_{x \in \Z}\min\left\{k>0\ :\ \lo([x-k,x+k],n)>n/2\right\} \leq c_{1}\right\} \\
& \subset & \left\{ \lo^*(n) > \frac{n}{2c_{1}}  \right\},
\end{eqnarray*}
we conclude with Theorem \ref{th1}.  
\begin{Rem} \label{rem1} We have seen that Theorem \ref{th2} implies Theorem \ref{th0}, moreover thanks to the result of \cite{GanShi} extended to Sinai's walk we know  that $\p\left[ \limsup_{n} \frac{ \lo^*\left(n\right)}{n}= \textrm{const} \in ]0,\infty]  \right]=1$, therefore there exists $c_{2}>0$ and $c_{3}>0$ such that for all $0 \leq \beta<1$
\begin{eqnarray}
\p\left[\limsup_{n}\left\{Y_{n,\beta}\leq c_{3}(1-\beta)^{-2},\ \lo^*(n)/n>c_{2} \right\}\right] =1. \label{a2}
\end{eqnarray}   
\end{Rem}

\section{Proof of Theorem \ref{th3} }

To prove this Theorem we use Remark \ref{rem1}. We begin with the following nice facts, define 
\begin{eqnarray}
& & \A_{n}= \max_{(x,y)\in \F_{n}^2} |x-y| \leq c_{4}, \\
& & \B_{n}= \max_{x} \lo\left(\left[x-c_{4}/2,x+c_{4}/2\right],n\right) >n-\lo^{*}(n),
\end{eqnarray}
recalling that $\lo^{*}(n)=\max_{x} \lo\left(x,n\right)$, $c_{4}$ is for the moment a free parameter that will be chosen latter.
\textit{Fact} 1 We have $\B_{n} \subseteq \A_{n} $, indeed it is easy to check that
\begin{eqnarray}
 \bigcap_{x \in \F_{n}}\left\{ \sum_{k=x-c_{4}/2}^{x+c_{4}/2}\lo\left(k,n\right) >n-\lo^{*}(n)\right\} \subseteq \A_{n}
\end{eqnarray}
moreover $\F_{n}\subset \Z$, so it is clear  that
\begin{eqnarray}
 \bigcap_{x \in \F_{n}}\left\{ \sum_{k=x-c_{4}/2}^{x+c_{4}/2}\lo\left(k,n\right) >n-\lo^{*}(n)\right\} \supseteq  \bigcap_{x \in \Z}\left\{ \sum_{k=x-c_{4}/2}^{x+c_{4}/2}\lo\left(k,n\right) >n-\lo^{*}(n)\right\} \equiv \B_{n}.
\end{eqnarray}
\textit{Fact} 2 Using \ref{a2} with $\beta=1-c_{2}$ we have
\begin{eqnarray}
\p\left[\limsup_{n} \left\{ Y_{n,(1-c_{2})} \leq c_{3}(c_{2})^{-2} \textrm{ and } \lo^*(n)/n \geq c_{2} \right\} \right]=1 \label{leeq}
\end{eqnarray}
Now, using the Definition of the ''$\liminf$'' and Fact 1 we get
\begin{eqnarray}
\p\left[\liminf _{n} \max_{(x,y)\in \F_{n}^2} |x-y| \leq c_{4} \right] \geq \p\left[\limsup_{n}  \B_{n}  \right].
\end{eqnarray}
It is clear that 
\begin{eqnarray}
\p\left[ \limsup_{n}  \B_{n}  \right] & \geq & \p\left[\limsup_{n} \left\{B_{n}\ , \frac{\lo^*(n)}{n} \geq c_{2} \right\}  \right], \end{eqnarray}
 moreover 
\begin{eqnarray}
\qquad   \left\{ B_{n}\ , \frac{\lo^*(n)}{n} \geq c_{2} \right\} \supseteq \left\{ \max_{x} \lo\left(\left[x-c_{4}/2,x+c_{4}/2\right],n\right) >n(1-c_{2}),\  \frac{\lo^*(n)}{n} \geq c_{2} \right\}.
\end{eqnarray}
Therefore choosing $c_{4}=c_{3}/(c_{2}^{2})$, we finally get that :
\begin{eqnarray}
\qquad \p\left[ \limsup_{n}  \B_{n}  \right] & \geq & \p\left[ \limsup_{n} \left\{ \max_{x} \lo\left(\left[x-c_{3}/c_{2}^{2},x+c_{3}/c_{2}^{2}\right],n \right) >n(1-c_{2}),\ \frac{\lo^*(n)}{n} \geq c_{2} \right\} \right] \\
&\equiv & \p\left[ \limsup_{n} \left\{Y_{n,(1-c_{2})} \leq c_{3}(c_{2})^{-2}  ,\ \frac{\lo^*(n)}{n} \geq c_{2} \right\} \right] \\
&=& 1
\end{eqnarray}
where the last equality comes from Fact 2. $\blacksquare$

\section{Conclusion remarks}

We have seen that using the method of \cite{Pierre2} and the Proposition 3.1 of  \cite{GanShi} we get easily annealed result for the concentration variable $Y_{n}$. We also point out that the result on the concentration variable implies both results on the maximum of the local time and on the favorite sites.

Here we only get the ''$\liminf$'' asymptotic of $Y_{n}$, what can we say about the ''$\limsup$'' ? We notice that if we have something like $\liminf \lo^{*}(n)\phi(n)/n=cte>0\ \p.a.s$  then $\limsup Y_{n}/\phi(n)=cte \in ]0+\infty],\ \p.a.s$ but is $\phi(n)$ the good asymptotic for the ''$\limsup$'' of $Y_{n}$ ? Notice that forthcoming work of Gantert shows that $\liminf \lo^{*}(n) \log \log \log n /n=cte>0\ \p.a.s$ and  forthcoming work of Z. Shi and O. Zindy implies that $\limsup Y_{n}/\log \log \log n=cte \in ]0+\infty],\ \p.a.s$.

Now, forgetting the hypothesis \ref{hyp1} and using the ones of \cite{GanShi} (originally introduced by \cite{KesKozSpi}) 
\begin{eqnarray}
-\infty < \E_Q\left[ \log
\frac{1-\alpha_0}{\alpha_0}\right]<0 , 
\end{eqnarray}
and that there is $0<\kappa <1$ such that 
\begin{eqnarray}
0< \E_Q\left[ \left(
\frac{1-\alpha_0}{\alpha_0}\right)^{\kappa}\right]=1.
\end{eqnarray}
Thanks to their work, it appears clearly that for small $\beta$ one can find $c_{1}\equiv c_{1}(\beta)>0$ such that  :
\begin{eqnarray}
\liminf Y_{n,\beta} \leq c_{1},\ \p.a.s \label{lasteq11}
\end{eqnarray}
a question that is maybe interesting is to understand how this $\beta$ depends on $\kappa$, for example, can we find $\kappa$ such that \ref{lasteq11} is true for $\beta=1/2$ ?
We could say that Sinai's walk is concentrated uniformly for $0<\beta<1$ whereas Kesten et al. walk is uniformly concentrated  for $0<\beta<\beta_{c}\equiv \beta_{c}(\kappa)$. What can we say about  $\beta_{c}$ ?

\appendix

\section{Basic results for birth and death processes}

\noindent For completeness we recall some results of
\cite{Chung} and \cite{Revesz} on inhomogeneous discrete time
birth and death processes.

\noindent Let $x,a$ and $b$ in $\Z$, assume $a<x<b$,
the two following lemmata can be found in \cite{Chung} (pages
73-76), the proof follows from the method of difference equations.

\begin{Lem} \label{3.7bb} Recalling \ref{3.7sec}, for all $\alpha$ we have
 \begin{eqnarray}
 & &
 \p^{\alpha}_x\left[T_a>T_b\right]=\frac{\sum_{i=a+1}^{x-1}\exp
\big(S_{i}-  S_{a}\big)  +1}{\sum_{i=a+1}^{b-1}\exp\big( S_{i}-
S_{a}  \big)+1 } \label{k1} , \\  & & \p^{\alpha}_x \left[T_a<T_b
\right]=\frac{\sum_{i=x+1}^{b-1} \exp \big( S_{i}- S_{b} \big)
+1}{\sum_{i=a+1}^{b-1}\exp \big( S_{i}- S_{b} \big) +1 }
\label{k2} .
\end{eqnarray}
\end{Lem}

\noindent Now we give some explicit expressions for the local
times that can be found in \cite{Revesz} (page 279)

\begin{Lem} \label{Lou} For all $\alpha$ and $i \in \Z$, we have, if $x>i$
\begin{eqnarray}
\Ea_i\left[
\lo(x,T_i)\right]=\frac{\alpha_i\p_{i+1}^{\alpha}\left[T_{x}<T_{i}\right]}{
\beta_{x}\p^{\alpha}_{x-1}\left[T_{x}>T_{i}\right]} ,
\label{2eq65}
\end{eqnarray}
if $x<i$
\begin{eqnarray}
\Ea_i\left[ \lo(x,T_i)\right]= \frac
{\beta_i\p_{i-1}^{\alpha}\left[T_{x}<T_{i}\right]}{
\alpha_{x}\p^{\alpha}_{x+1}\left[T_{x}>T_{i}\right]} .
\label{2eq68}
\end{eqnarray}
\end{Lem}


\section{Proof of the good properties for the environment \label{sec4}}

Here we give the main ideas for the proof of the Proposition \ref{profondab}, we begin with some

\subsection{Elementary results \label{ER1} for sum of i.i.d. random variables}

We will always work on the right hand side of the origin, that
means with $(S_m,m\in \N)$, by symmetry we obtain the same result
for $m \in \Z_-$.

 \noindent \\ We introduce the following stopping times, for
$a>0$,
\begin{eqnarray}
& & V^+_a\equiv V^+_a(S_j,j\in \N) = \left\{\begin{array}{l}  \inf  \{m \in \N^*,\ S_m \geq a \} , \label{def1}  \\
 + \infty \textrm{, if such a } m \textrm{ does not exist.} \end{array} \right. \label{4.3} \\
& & V^-_a\equiv V^-_a(S_j,j\in \N)= \left\{\begin{array}{l}  \inf  \{m \in \N^*,\ S_m \leq -a \} ,   \\
 + \infty \textrm{, if such a } m \textrm{ does not exist.}
\end{array} \right.
\end{eqnarray}
The following lemma is an immediate consequence of the Wald
equality (see \cite{Neveu})
\begin{Lem} \label{lem101b}  Assume \ref{hyp1bb}, \ref{hyp0} and \ref{hyp4}, let $a>0$, $d>0$ we have
\begin{eqnarray}
& & Q\left[V_{a}^{-} < V_{d}^{+} \right] \leq
\frac{d+\Ie}{d+a+\Ie} , \label{lem101eq2b} \\
& & Q\left[V_{a}^{-} > V_{d}^{+} \right] \leq
\frac{a+\Ie}{d+a+\Ie} , \label{lem101eq3b}
\end{eqnarray}
with $ \Ie \equiv \log ((1-\eto)(\eto)^{-1})$.
\end{Lem}

\noindent The following lemma is a basic fact for sums of i.i.d. random variables 

\begin{Lem} \label{maldita} Assume \ref{hyp1bb}, \ref{hyp0} and \ref{hyp4}
hold, there exists $b \equiv b(Q)>0$  such that for all $r>0$
\begin{eqnarray}
Q\left[V^-_0>r \right] \leq \frac{b}{\sqrt{r}}
\label{8lem29}.
\end{eqnarray}
\end{Lem}

\subsection{Proof of Proposition \ref{profondab} }

It is in this part where the $(1-\beta)^{-2}$ dependance occuring in Theorem \ref{th2}  will become clear. The main difficulty is to get an upper bound for the expectation $\E_{Q}\left[\Ea_{\tmo}\left[ \lo(\tilde{\Theta}(n,\beta),T_{\tmo})\right]\right] $.

\subsubsection{Preliminaries  \label{9.1.2}} \noindent

\noindent
By linearity of the expectation we have :
\begin{eqnarray}
\qquad \Ea_{\tmo}\left[ \lo(\tilde{\Theta}(n,\beta),T_{\tmo})\right] &
 \equiv & \sum_{j=\tmo+\tcu}^{M_n}\Ea_{\tmo}\left[ \lo(j,T_{\tmo})\right] +
\sum_{j=M_n'}^{\tmo-\tcu}\Ea_{\tmo}\left[ \lo(j,T_{\tmo})\right] +1 \label{eq23},
\end{eqnarray} 
recall that $\tcu=c_{3}(1-\beta)^{-2}$ with $c_{3}>0$ and $0\leq \beta<1$. Now using
Lemma \ref{3.7bb} and hypothesis \ref{hyp4} we easily get the
following lemma
\begin{Lem} \label{9eq13}  Assume \ref{hyp4}, for all  ${M_n}' \leq k \leq {M_n}$, $k\neq \tmo$
\begin{eqnarray}
\frac{\eto}{1-\eto}\frac{1}{e^{S_k-S_{\tmo}}} \leq
\Ea_{\tmo}\left[ \lo(k,T_{\tmo})\right] \leq
\frac{1}{\eto}\frac{1} {e^{S_k-S_{\tmo}} } ,
\end{eqnarray}
with a $Q$ probability equal to one.
\end{Lem}\noindent 
The following lemma is easy to prove :
\begin{Lem} \label{3lem62} For all $n>3$, with a $Q$ probability equal to one we have
\begin{eqnarray}
& & \sum_{j=\tmo+\tcu}^{M_n}\frac{1}{e^{S_{j}-S_{\tmo}}} \leq
\sum_{i=1}^{N_{n}+1} \frac{1}{e^{a(i-1)}}\sum_{j=m_n+\tcu}^{M_n}
\un_{S_{j}-S_{\tmo} \in [a(i-1),ai[} ,\label{3eqlem1.209} \\
& & \sum_{j=M_n'}^{m_n-\tcu} \frac{1}{e^{S_{j}-S_{\tmo}}} \leq
\sum_{i=1}^{N_{n}+1}
\frac{1}{e^{a(i-1)}}\sum_{j=M_n'}^{m_n-\tcu} \un_{S_{j}-S_{\tmo} \in
[a(i-1),ai[} ,\label{3eqlem1.210}
\end{eqnarray}
where $a=\frac{\Ie}{4}$, $N_{n}=[(\Gamma_n+\Ie)/a]$, recall that  $\un$ is the indicator function.
\end{Lem}

\noindent \\ Using \ref{eq23}, Lemma \ref{9eq13} and \ref{3lem62}, we have for
all $n>3$
\begin{eqnarray}
 \E_Q\left[ \Ea_{\tmo}\left[ \lo(\tilde{\Theta}(n,\beta),T_{\tmo})\right] \right]
\leq  1 &+& \frac{1}{\eta_{0}}\sum_{i=1}^{N_{n}+1} \frac{1}{e^{a(i-1)}}
\E_Q\left[\sum_{j=m_n+\tcu}^{M_n} \un_{S_{j}-S_{\tmo} \in
[a(i-1),ai[} \right] \label{9.161} \\ &+& \frac{1}{\eta_{0}}
\sum_{i=1}^{N_{n}+1} \frac{1}{e^{a(i-1)}}\E_Q\left[
\sum_{j=M_n'}^{m_n-\tcu} \un_{S_{j}-S_{\tmo} \in [a(i-1),ai[}\right]
. \nonumber
\end{eqnarray}

\noindent The next step for the proof is to show that the two expectations $\E_Q[...]$ on
the right hand side of \ref{9.161} are bounded by a constant
depending only on the distribution $Q$ times a polynomial in $i$ times $1/\sqrt{\tcu}$:

\begin{Lem} \label{3lem} There exits a constant $c\equiv c(Q)$ such that for all $n$ large enough :
\begin{eqnarray}
& & \E_Q\left[\sum_{j=m_n+\tcu}^{M_n} \un_{S_{j}-S_{\tmo} \in
[a(i-1),ai[} \right] \leq \frac{c \times i^3}{\sqrt{\tcu}} \label{9.161bb}, \\
& & \E_Q\left[ \sum_{j=M_n'}^{m_n-\tcu} \un_{S_{j}-S_{\tmo} \in
[a(i-1),ai[}\right] \leq \frac{c \times i^3}{\sqrt{\tcu}}.  \label{9.161bbb}
\end{eqnarray}
\end{Lem}

\subsubsection{Proof of Lemma \ref{3lem}}

\begin{Rem}
We give some details of the proof of Lemma \ref{3lem}
mainly because it helps to understand the appearance of the $(1-\beta)^{-2}$ in Theorem \ref{th2}, moreover it is based on a very nice cancellation that occurs
between two $\Gamma_{n} \equiv \log n+ \log_{2}n$, see formulas \ref{5.29} and \ref{5.69}.
Similar cancellation is already present in \cite{Kesten2}.
\end{Rem}

\noindent Let us define the following stopping times, let $i>1$ :
\begin{eqnarray*}
& & u_0=0, \\
& & u_1 \equiv V_0^-= \inf\{m>0,\ S_m<0\}, \\
& & u_i = \inf\{m>u_{i-1},\ S_m<S_{u_{i-1}}\}.
\end{eqnarray*}

\noindent \\ The following lemma give a way to characterize the
 point $\tmo$, it is inspired by the work of \cite{Kesten2} and is
 just inspection

\begin{Lem} \label{8eq26b} Let $n>3$ and $\gamma>0$, recall $\Gamma_n =\log n +\gamma \log_2 n$, assume $\tmo > 0 $, for all  $l\in
\N^*$ we have
\begin{eqnarray}
 \tmo=u_l &  \Rightarrow & \left\{
\begin{array}{l}
\bigcap_{i=0}^{l-1}\left\{ \max_{ u_i \leq j \leq u_{i+1} } (S_i)-S_{u_{i}} < \Gamma_n\right\} \textrm{ and } \\
\max_{ u_l \leq j \leq u_{l+1} } (S_i)-S_{u_{i}} \geq \Gamma_n
 \textrm{ and } \\
{M_n}=V^+_{\Gamma_n,l}
\end{array} \right. \label{2eq121b}
\end{eqnarray}
where
\begin{eqnarray}
V^+_{z,l}\equiv V^+_{z,l}\left(S_j,j\geq 1\right) =
\inf\left(m>u_l,\ S_m-S_{u_l} \geq z  \right).
\end{eqnarray}
\end{Lem}

\noindent A similar characterization of $\tmo$ if $\tmo\leq 0$ can
be done (the case $\tmo=0$ is trivial).
We will only prove \ref{9.161bb}, we get \ref{9.161bbb} symmetrically 
moreover we assume that $m_n>0$, computations are similar for the
case $m_n \leq 0$. Thinking on the basic definition of the expectation, we 
 need an upper bound for the probability :
\begin{eqnarray*}
Q\left[\sum_{j=m_n+\tcu}^{M_n} \un_{S_{j}-S_{\tmo} \in
[a(i-1),ai[}=k. \right]
\end{eqnarray*}
First we make a partition over the values of $m_n$ and then we
use Lemma \ref{8eq26b}, we get  :
\begin{eqnarray} 
\quad Q\left[\sum_{j=m_n+\tcu}^{M_n} \un_{S_{j}-S_{\tmo} \in
[a(i-1),ai[}=k. \right] & \equiv &  \sum_{l \geq 0}
Q\left[\sum_{j=m_n+\tcu}^{M_n} \un_{S_{j}-S_{\tmo} \in
[a(i-1),ai[}=k,\ \tmo=u_l \right] \nonumber \\
 & \leq & \sum_{l \geq 0} Q\left[ \A_{\Gamma_n,l}^+,\max_{ u_l \leq j \leq u_{l+1} } (S_j)-S_{u_{l}}
\geq \Gamma_n,\A_{\Gamma_n,l}^-  \right] \label{5.64}
\end{eqnarray}
where
\begin{eqnarray*}
& & \A_{\Gamma_n,l}^+=\sum_{s=u_l+\tcu}^{V^+_{\Gamma_n,l}}\un_{\{S_j-S_{u_l} \in
[a(i-1),ai[\}}=k,\ \\
& & \A_{\Gamma_n,l}^- =\bigcap_{r=0}^{l-1}\left\{ \max_{ u_r \leq j \leq
u_{r+1} } (S_r)-S_{u_{r}} < \Gamma_n\right\},\ \A_0^-=\Omega_1.
\end{eqnarray*}
for all $l\geq 0$. By the strong Markov property we have :
\begin{eqnarray}
 Q\left[ \A_{\Gamma_n,l}^+,\ \max_{ u_l \leq j \leq u_{l+1} } (S_j)-S_{u_{l}}
\geq \Gamma_n,\ \A_{\Gamma_n,l}^-  \right]  \leq  Q\left[ \A_{\Gamma_n,0}^+,\
V_0^->V^+_{\Gamma_n} \right] Q\left[\A_{\Gamma_n,l}^- \right]. \label{5.63}
\end{eqnarray}
The strong Markov property gives also that the sequence $(\max_{
u_r \leq j \leq u_{r+1} } (S_r)-S_{u_{r}} < \Gamma_n ,r \geq 1)$
 is i.i.d., therefore :
\begin{eqnarray}
Q\left[\A_{\Gamma_n,l}^-
\right] \leq \left(Q\left[V_0^-<V^+_{\Gamma_n}\right]\right)^{l-1}. \label{5.60}
\end{eqnarray}
We notice that $Q\left[ \A_{\Gamma_n,0}^+,\ V_0^->V^+_{\Gamma_n} \right]$
 does not depend on $l$, therefore, using \ref{5.64}, \ref{5.63}
 and \ref{5.60} we get :
\begin{eqnarray}
\quad Q\left[\sum_{j=m_n+\tcu}^{M_n} \un_{S_{j}-S_{\tmo} \in
[a(i-1),ai[}=k. \right] & \leq & (1+(Q\left[V_0^-\geq
V^+_{\Gamma_n}\right])^{-1})Q\left[ \A_{\Gamma_n,0}^+,\ V_0^->V^+_{\Gamma_n}
\right] \label{5.29}.
\end{eqnarray}
Using the Markov property we obtain that  
\begin{eqnarray}
& & Q\left[ \A_{\Gamma_n,0}^+,\ V_0^->V^+_{\Gamma_n} \right] \leq Q\left[V^{-}_{0}>\tcu\right]\max_{ 0 \leq x \leq \tcu/\Ie} \left\{ Q_{x}\left[ \A_{\Gamma_n,0}^+,\ V_0^->V^+_{\Gamma_n} \right] \right\}.
\end{eqnarray}
$\Ie$ is given just after \ref{lem101eq3b}. To get an upper bound for $Q_{x}\left[ \A_{\Gamma_n,0}^+,\ V_0^->V^+_{\Gamma_n}
\right]$, we introduce the following sequence of stopping
times, let $k>0$ :
\begin{eqnarray*}
& & H_{ia,0}=0, \\
& & H_{ia,k}=\inf\{m>H_{ia,k-1},\ S_m \in [(i-1)a,ia[\}.
\end{eqnarray*}
Making a partition over the values of $H_{ia,k}$ and using the Markov property we get:
\begin{eqnarray}
& & Q_{x}\left[ \A_{\Gamma_n,0}^+,\ V_0^->V^+_{\Gamma_n} \right] \nonumber \\
& \leq & \sum_{w \geq 0}\int_{(i-1)a}^{ia} Q_{x}\left[H_{ia,k}=w,S_w \in
dy,\bigcap_{s=0}^{w}\{S_s>0\},\bigcap_{s=w+1}^{\inf \{l>w,S_l
\geq \Gamma_n-x\}}\{S_s>0\} \right]  \nonumber \\
& \leq &  Q_{x}\left[H_{ia,k}<V^-_0\right] \max_{ (i-1)a \leq y \leq
ia } \left\{Q_y\left[V_{\Gamma_n-y}^{+}<V^-_{y} \right]\right\} \nonumber \\
& \equiv & Q_{x}\left[H_{ia,k}<V^-_0\right] Q_{ia}\left[V_{\Gamma_n-ia}^{+}<V^-_{ia}
\right] \label{5.69}
\end{eqnarray}
To finish we need an upper bound for
$Q_{x}\left[H_{ia,k}<V^-_0\right]$, we do not want to give details of
the computations for this because it is not difficult, however the
reader can find these details in \cite{Pierreth} pages 142-145. We
have for all $i>1$:
\begin{eqnarray}
\qquad  Q_{x}\left[H_{ia,k}<V^-_0\right]  &\leq&
Q_{x}\left[V^-_0>V^+_{(i-1)a}\right]
\left(1-Q\left[\epsilon_0<-\frac{\Ie}{2}\right]Q_{(i-1)a-\frac{\Ie}{4}}\left[V^+_{(i-1)a}
\geq V^-_0\right] \right)^{k-1} \  \nonumber  \\ 
& \leq & 
\left(1-Q\left[\epsilon_0<-\frac{\Ie}{2}\right]Q_{(i-1)a-\frac{\Ie}{4}}\left[V^+_{(i-1)a}
\geq V^-_0\right] \right)^{k-1},
\label{5.70} 
\end{eqnarray}
and in the same way
\begin{eqnarray}
Q_{x}\left[H_{a,k}<V^-_0\right]  \leq 
\left(1-Q\left[\epsilon_0<-\frac{\Ie}{4}\right] \right)^{k-1} \ .  \label{5.68}
\end{eqnarray}
So using \ref{5.29}-\ref{5.68}, Lemmata \ref{lem101b} and \ref{maldita} one can find a constant $c \equiv c(Q)$ that
depends only on the distribution $Q$ such that for all $i\geq 0$ : 
\begin{eqnarray*}
& & \E_Q\left[\sum_{j=m_n+\tcu}^{M_n} \un_{S_{j}-S_{\tmo} \in
[a(i-1),ai[} \right]  \equiv  \sum_{k=1}^{+ \infty  } k
Q\left[\sum_{j=m_n+\tcu}^{M_n} \un_{S_{j}-S_{\tmo} \in
[a(i-1),ai[}=k \right]    \leq \frac{c \times  i^3}{\sqrt{\tcu}},
\end{eqnarray*}
which provide \ref{9.161bb}. $\blacksquare$
\noindent \\
Using both \ref{9.161} and Lemma \ref{3lem} we get that there exists $c_{0}\equiv c_{0}(Q)$ such that
\begin{eqnarray}
\qquad \E_{Q} \left[\Ea_{\tmo}\left[ \lo(\tilde{\Theta}(n,\beta),T_{\tmo})\right] \right] &
 \leq & \frac{c_{0}}{\sqrt{\tcu}}.
\end{eqnarray}
Now, using the elementary Markov inequality, we get :
\begin{eqnarray*}
& & Q\left[ \Ea_{\tmo}\left[ \lo(\tilde{\Theta}(n,\beta),T_{\tmo})\right] \leq  2 \frac{c_{0}}{\sqrt{\tcu}} \right] \\
& \equiv & 1-Q\left[ \Ea_{\tmo}\left[ \lo(\tilde{\Theta}(n,\beta),T_{\tmo})\right] >2  \frac{c_{0}}{\sqrt{\tcu}} \right]  \\
& \geq & 1/2 .
\end{eqnarray*}

{\small \bibliography{article}}

\begin{thebibliography}{16}
\expandafter\ifx\csname natexlab\endcsname\relax\def\natexlab#1{#1}\fi
\expandafter\ifx\csname url\endcsname\relax
  \def\url#1{{\tt #1}}\fi

\bibitem[Solomon(1975)]{Solomon}
F.~Solomon.
\newblock Random walks in random environment.
\newblock {\em Ann. Probab.}, \textbf{3}\penalty0 (1):\penalty0 \ 1--31, 1975.

\bibitem[Sinai(1982)]{Sinai}
Ya.~G. Sinai.
\newblock The limit behaviour of a one-dimensional random walk in a random
  medium.
\newblock {\em Theory Probab. Appl.}, \textbf{27}\penalty0 (2):\penalty0 \
  256--268, 1982.

\bibitem[Deheuvels and R\'ev\'esz(1986)]{Deh&Revesz}
P.~Deheuvels and P.~R\'ev\'esz.
\newblock Simple random walk on the line in random environment.
\newblock {\em Probab. Theory Relat. Fields}, \textbf{72}:\penalty0 \ 215--230,
  1986.

\bibitem[Hu and Shi(1998)]{HuShi2}
Y.~Hu and Z.~Shi.
\newblock The limits of {S}inai's simple random walk in random environment.
\newblock {\em Ann. Probab.}, \textbf{26}\penalty0 (4):\penalty0 \ 1477--1521,
  1998.

\bibitem[R\'ev\'esz(1989)]{Revesz}
P.~R\'ev\'esz.
\newblock {\em Random walk in random and non-random environments}.
\newblock World Scientific, 1989.

\bibitem[R\'ev\'esz(1988)]{Revesz1}
P.~R\'ev\'esz.
\newblock In random environment the local time can can be very big.
\newblock {\em Ast\'erisque}, pages \ 157--158,321--339, 1988.

\bibitem[Shi(1998)]{Shi}
Z.~Shi.
\newblock A local time curiosity in random environment.
\newblock {\em Stoch. Proc. Appl.}, \textbf{76}\penalty0 (2):\penalty0 \
  231--250, 1998.

\bibitem[Hu and Shi(2000)]{HuShi0}
Y.~Hu and Z.~Shi.
\newblock The problem of the most visited site in random environment.
\newblock {\em Probab. Theory Relat. Fields}, \textbf{116}\penalty0
  (2):\penalty0 \ 273--302, 2000.

\bibitem[Shi(2001)]{Shi1}
Z.~Shi.
\newblock {S}inai's walk via stochastic calculus.
\newblock {\em Panoramas et Synth\`eses}, \textbf{12}:\penalty0 \ 53--74, 2001.

\bibitem[Andreoletti(2005)]{Pierre2}
P.~Andreoletti.
\newblock On the concentration of {S}inai's walk.
\newblock {\em to appear in Stoch. Proc. Appl.}, 2005.

\bibitem[Gantert and Shi(2002)]{GanShi}
N.~Gantert and Z.~Shi.
\newblock Many visits to a single site by a transient random walk in random
  environment.
\newblock {\em Stoch. Proc. Appl.}, 99:\penalty0 \ 159--176, 2002.

\bibitem[Kesten(1986)]{Kesten2}
H.~Kesten.
\newblock The limit distribution of {S}inai's random walk in random
  environment.
\newblock {\em Physica}, \textbf{138A}:\penalty0 \ 299--309, 1986.

\bibitem[Kesten et~al.(1975)Kesten, Kozlov, and Spitzer]{KesKozSpi}
H.~Kesten, M.V. Kozlov, and F.~Spitzer.
\newblock A limit law for random walk in a random environment.
\newblock {\em Comp. Math.}, \textbf{30}:\penalty0 \ 145--168, 1975.

\bibitem[Chung(1967)]{Chung}
K.~L. Chung.
\newblock {\em Markov Chains}.
\newblock Springer-Verlag, 1967.

\bibitem[Neveu(1972)]{Neveu}
J.~Neveu.
\newblock {\em Martinguales \`a temps discret}.
\newblock Masson et Cie, 1972.

\bibitem[Andreoletti(2003)]{Pierreth}
P.~Andreoletti.
\newblock {\em Localisation et Concentration de la marche de {S}inai}.
\newblock PhD thesis, Universit\'e d'Aix-Marseille II, France, 2003.

\end{thebibliography}

\vspace{1cm} \noindent
\begin{tabular}{l}
Laboratoire Analyse-Topologie-Probabilit\'es - C.N.R.S. UMR 6632  \\
Centre de math\'ematiques et d'informatique \\
Universit\'e de Provence, \\
39 rue F. Joliot-Curie, \\
13453 Marseille cedex 13 \\
France
\end{tabular}

\end{document}